\documentclass[12pt,fleqn]{article}
\usepackage{amssymb}
\usepackage{times}
\usepackage{amsmath}
\usepackage{times}
\usepackage{url}
\usepackage{graphicx}
\begin{document}

\title{Fast Multipole Method based filtering of non-uniformly sampled data}
\date{October 30, 2003}
\author{Nail A. Gumerov, Ramani Duraiswami \thanks{We gratefully acknowledge support of NSF awards 0219681 and  0086075.} \\ Perceptual Interfaces \& Reality Laboratory, \\ UMIACS, University of Maryland, College Park}
\maketitle
\abstract{Non-uniform fast Fourier Transform (NUFFT) and inverse NUFFT (INUFFT) 
algorithms, based on the Fast Multipole Method (FMM) are developed and tested. Our 
algorithms are based on a novel factorization of the FFT kernel, and are implemented
with attention to data structures and error analysis. \\ \\
{\em Note: This unpublished manuscript was available on our web pages and has been referred to by others in the literature. To provide a proper archival reference we are placing it on arXiv.}}
\section{Introduction}

Most signal processing involves the manipulation of various transforms of
discretely sampled data, especially of Fourier transforms. Direct
computation of a transform of $N$ data samples requires $ON^{2})$
operations. For band-limited data that is sampled at equispaced locations,
the rediscovery in 1965 by  \cite{CooleyTukey65} of an $O(N\log N)$
algorithm, the fast Fourier transform, lead to a revolution. It is difficult
to overstate its impact on science and industry.

Unfortunately, many applications require transforms between nonuniformly
sampled data, and recently there has been much interest in this area,
especially for various inverse scattering calculations and as a component of
the computation of other transforms, e.g., on the sphere. The goal is to
achieve fast ($o\left( N^{2}\right) )$ bandwidth-preserving transforms of
data. Due to reasons of space we do not discuss these here, and refer the
reader to a recent review of NUFFTs \cite{Potts+01}. We focus instead on
extensions of the work of Dutt \&Rokhlin \cite{DuttRokhlin95}, who used the
Fast Multipole Method (FMM) to achieve a NUFFT\ of asymptotic complexity $%
O(N(\log N-\log \epsilon ))$, where $\epsilon $ is the specified error.

Here, we too apply the FMM to obtain NUFFTs and INUFFTs. The major
difference is in the factorization we used for the kernel function, and our
efficient implementation using data structures, supported by tests for $N$ $%
\simeq 10^{6}$. We show that for machine precision, the complexity of the
FMM-based NUFFT and INUFFT is $O(N)$. We present theoretical and
computational analyses of error bounds and optimal settings for the FMM. 

\section{Formulation \& Solution}

We consider the following two problems: i) a 2$\pi $-periodic complex-valued
band-limited function $f(x)$ is sampled at $N$ equispaced points, $%
x_{k}=2\pi k/N,$ $k=0,...,N-1,$ so $f_{k}=f(x_{k})$; \ find $g_{j}=f(y_{j}),$
$j=0,...,M-1$, where $y_{j}\in \left[ 0,2\pi \right) $ are non-equispaced
points with specified accuracy $\epsilon ;$ ii) the same function is given
at $M$ non-uniform samples $g_{j},$ $j=0,...,M-1$ on a non-uniform grid $%
\left\{ y_{j}\right\} $; determine $f_{k}$, $k=0,...,N-1$, \ on a uniform
grid $\left\{ x_{k}\right\} .$ $N$ determines the \emph{bandwidth} of $f$,
since 
\begin{equation}
f(x)=\sum_{n=0}^{N-1}c_{n}e^{inx},\quad c_{n}=\frac{1}{N}%
\sum_{k=0}^{N-1}f_{k}e^{-inx_{k}},  \label{p1}
\end{equation}%
where $c_{n}$ are the Fourier coefficients. The forward DFT is the transform 
$\left\{ f_{k}\right\} \rightarrow \left\{ c_{n}\right\} ,$ and the IDFT is
transform $\left\{ c_{n}\right\} \rightarrow \left\{ f_{k}\right\} $ . \ For 
$N=2^{L},$ the transforms $\left\{ f_{k}\right\} \rightleftarrows \left\{
c_{n}\right\} $ can be done in $O(N\log N)$ operations using the FFT.
Therefore, if interpolation problems $\left\{ f_{k}\right\} \rightleftarrows
\left\{ g_{j}\right\} $ can be solved with the same efficiency, NUFFT and
INUFFT, which are transforms $\left\{ c_{n}\right\} \rightleftarrows $ $%
\left\{ g_{j}\right\} $ will be available. Eq. (\ref{p1}) can be written as 
\begin{equation}
\left\{ g_{j}\right\} =\left\{ K_{jk}\right\} \left\{ f_{k}\right\} ,\quad
K_{jk}=\frac{1}{N}\sum_{n=0}^{N-1}e^{-inx_{k}}e^{iny_{j}},  \label{p2}
\end{equation}%
where $j=0,...,M-1;$ $k=0,...,N-1.$ It relates the vectors $\left\{
g_{j}\right\} $ and $\left\{ f_{k}\right\} $ via the sampled kernel matrix $%
\left\{ K_{jk}\right\} $. The sum (\ref{p2}) is a geometric progression, and
can be simplified as 
\begin{subequations}
\begin{gather}
\!\!\!\!\!\!\!\!\!\!K_{jk}=\frac{e^{iNy_{j}}-1}{e^{i\left( y-x_{k}\right) }-1%
}=F_{j}G\left( y_{j}-x_{k}\right) ,  \label{p4} \\
\!\!\!\!\!\!\!\!\!\!F_{j}=\frac{e^{iNy_{j}}-1}{N},G\left( t\right) =\frac{1}{%
e^{it}-1}=-\frac{1+i\cot \left( t/2\right) }{2}.  \label{p5}
\end{gather}%
This decomposition was derived in \cite{DuttRokhlin95}. It can be thought of
as a modulation of a fast oscillating function $F(y)$ by a slowly changing
amplitude function $G\left( t\right) $ amplitude. Evaluation of $F(y_{j})$
costs $O(M)$, while evaluation of $G\left( t\right) $ can be done
approximately by fast methods, such as the FMM. The transform $\left\{
g_{j}\right\} \rightarrow \left\{ f_{k}\right\} $ can be reduced similarly.
It leads to 
\end{subequations}
\begin{equation}
\left\{ f_{k}\right\} =\left\{ L_{kj}\right\} \left\{ g_{j}\right\}
,~~L_{kj}=C_{k}G\left( x_{k}-y_{j}\right) D_{j},  \label{p6}
\end{equation}%
where $j=0,...,M-1;$ $k=0,...,N-1$. The kernel $L$ is decomposed as in \cite%
{DuttRokhlin95}, where $G\left( t\right) $ is given by Eq. (\ref{p5}); and $%
C_{k}$ and $D_{j}$ are coefficients depending on $\left\{ x_{k}\right\} $
and $\left\{ y_{j}\right\} $ 
\begin{equation}
\!\!\!\!\!\!C_{k}=\left( -1\right) ^{k}\!\prod_{j=0}^{N-1}\sin \frac{%
x_{k}-y_{j}}{2},\quad D_{j}\!=\!\frac{2ie^{-iNy_{j}/2}}{\prod_{\substack{ %
k=0,  \\ k\neq j}}^{N-1}\sin \frac{y_{j}-y_{k}}{2}}.  \label{p8}
\end{equation}%
These coefficients can be preccomputed using the FMM. So both the forward
and inverse problems can be reduced to multiplication of a matrix generated
by kernel $G\left( t\right) $ by some input vector. This is the problem
which is in the focus of the present paper. 

\subsection{Multilevel FMM}

For computation of the matrix-vector product $\left\{ G\left(
y_{j}-x_{k}\right) \right\} \left\{ f_{k}\right\} $ we use the Multilevel
FMM (MLFMM), a description of which can be found elsewhere \cite%
{Greengard88, GDB03}\textbf{. }Note that $\left\{ x_{k}\right\} $ or $%
\left\{ y_{j}\right\} $ need not be equispaced for the FMM. Both (\ref{p2})
and (\ref{p6}) use the same algorithm with source and target sets $\left\{
x_{k}\right\} $ \& $\left\{ y_{j}\right\} $ exchanged.

\begin{figure}[htb]
\includegraphics[width=0.9\textwidth,trim=0in 4in 0in 0in]{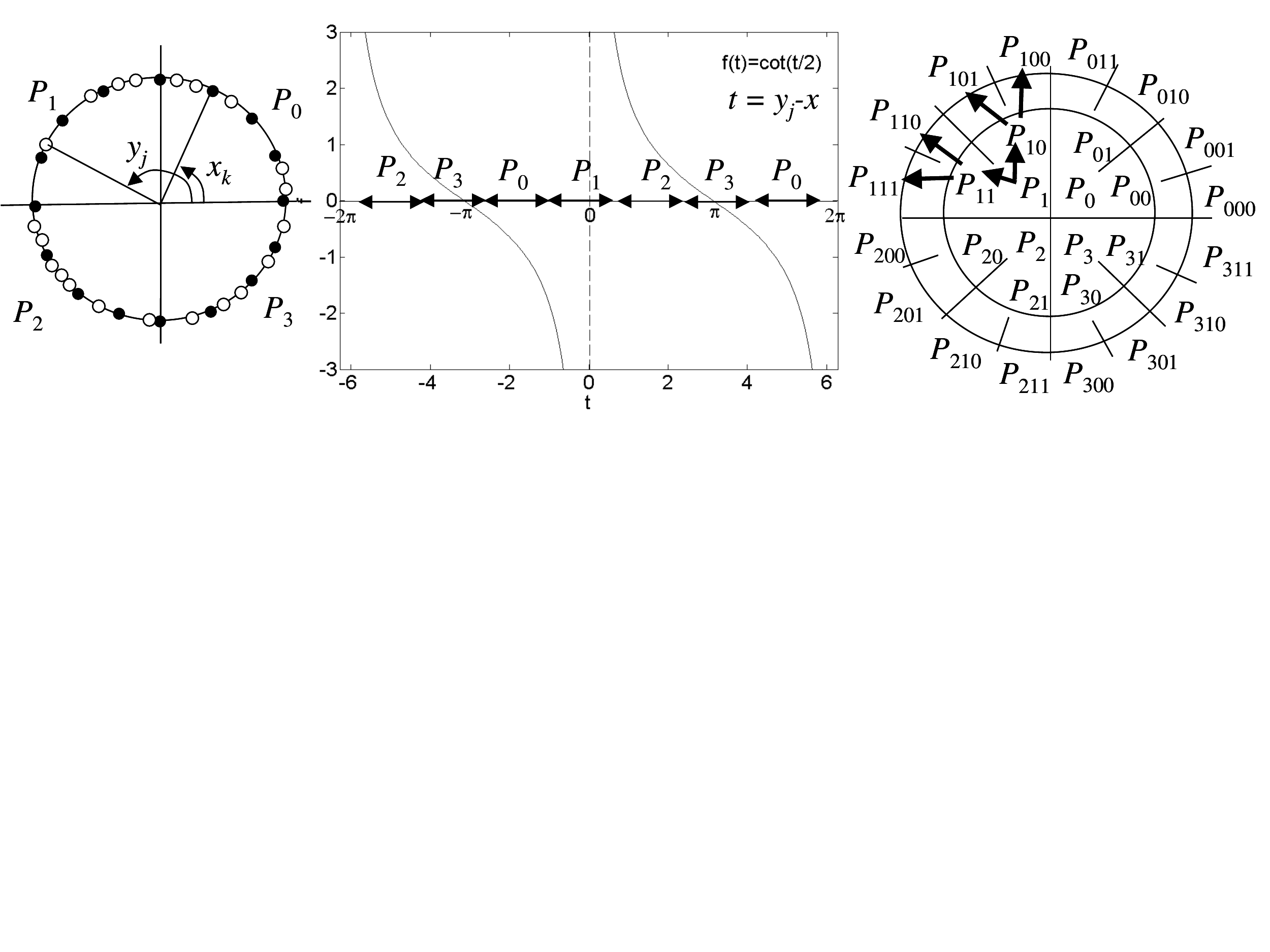}
\caption{On the left locations of the
	source, $\left\{ x_{k}\right\} ,$ and target, $\left\{ y_{j}\right\} ,$
	points on the unit circle are shown. The circle is subdivided at level 2
	into 4 arcs. In the center the behavior of the kernel function $\cot (t/2)$
	is shown. At the right the hierarchical space partitioning employed in the
	MLFMM is illustrated. Arrows show the hierarchical parent-children relations
	in the binary tree. \label{Fig.1}}
\end{figure}

\textbf{Space Partitioning:} The kernel $\cot \frac{y_{j}-x_{k}}{2}$ is a $%
2\pi $-periodic function. We can thus make transforms $\widetilde{x}%
_{k}=x_{k}+2\pi n,$ $n=0,\pm 1,...$, which do not change the function, and
keep $y_{j}-$ $\widetilde{x}_{k}$ in $-\pi \leqslant y_{j}-$ $\widetilde{x}%
_{k}\leqslant \pi ;,$ i.e., $y_{j}$ and $x_{k}$ are points on a unit circle.
Then $x_{k}$ and $\widetilde{x}_{k}$ are identical (see Fig. \ref{Fig.1}). 

In the case of circular topology, hierarchical space partitioning with a
binary tree with $l_{\max }$ levels$,$ leads to a division of the circle by $%
2^{l}$ equal arcs. In this hierarchy, parent-children relations are defined
as in a binary tree, but the concept of neighbors and neighborhoods is
modified. At level 2 we subdivide $\left[ 0,2\pi \right) $ into 4 pieces: $%
P_{0}:\left[ 0,\pi /2\right) ,$ $P_{1}:\left[ \pi /2,\pi \right) ,$ $P_{2}:%
\left[ \pi ,3\pi /2\right) ,$ $P_{3}:\left[ 3\pi /2,2\pi \right) ,$ where $%
P_{2}$ and $P_{0}$ are neighbors of $P_{3}$ and they together form a
neighborhood of $P_{3},$ which we denote as $\Omega _{1}(P_{3}).$ We also
denote $P_{1}$ as $\Omega _{2}(P_{3}),$ i.e. the piece outside the
neighborhood of $P_{3}$.

\subsection{Kernel Expansions: }
Level $l=2$ requires special consideration. Here
we separate the sum over all the sources as 
\begin{equation}
\!\!\!\!\!h\left( y_{j}\right) =\!\!\sum_{x_{k}\in \Omega
_{1}(P_{n})}f_{k}\cot \frac{y_{j}-x_{k}}{2}+\!\!\sum_{x_{k}\in \Omega
_{2}(P_{n})}f_{k}\cot \frac{y_{j}-x_{k}}{2},  \label{p9}
\end{equation}%
$y_{j}\in P_{n};$ $n=0,...,3.$ The first sum can be computed using: 
\begin{equation}
\cot \frac{t}{2}=\frac{2}{t}-2\sum_{m=1}^{\infty }\frac{\left\vert
B_{2m}\right\vert }{\left( 2m\right) !}t^{2m-1},\quad \left\vert
t\right\vert <2\pi ,  \label{10}
\end{equation}%
where $B_{m}$ are the Bernoulli numbers \cite{AS64}. When applied here, we
note that for $x_{k}\in \Omega _{1}(P_{n})$ we have $-\pi \leqslant y_{j}-%
\widetilde{x}_{k}\leqslant \pi ,$ so for $t$ $=y_{j}-\widetilde{x}_{k}$ we
have $\left\vert t\right\vert \leqslant \pi .$ This provides fast
convergence in Eq. (\ref{10}), which otherwise blows up as $%
y_{j}-x_{k}\rightarrow \pm 2\pi .$

For the second sum, we introduce $t=y_{j}-\widehat{x}_{k}$, where $\widehat{x%
}_{k}=x_{k}\pm \pi $ is a point opposite $x_{k}$. We have 
\begin{equation}
\!\!\cot \frac{y_{j}-x_{k}}{2}\!=\!-\tan \frac{t}{2}\!=\!-2\sum_{m=1}^{%
\infty }\frac{\left( 2^{2m}-1\right) \left\vert B_{2m}\right\vert }{\left(
2m\right) !}t^{2m-1},  \label{11}
\end{equation}%
$\left\vert t\right\vert <\pi .$ Note that as $t\rightarrow \pm \pi $ the
series blows up. However, since $x_{k}\in \Omega _{2}(P_{n})$ we have $%
\widehat{x}_{k}\in P_{n}.$ Due to both $y_{j}$ and $\widehat{x}_{k}$ belong
to the same $P_{n}$, whose length is $\pi /2,$ we have $-\pi /2\leqslant
y_{j}-\widehat{x}_{k}\leqslant \pi /2,$ or $\left\vert t\right\vert
\leqslant \pi /2.$ This provides a fast convergence of the sum.

There are thus two parts to the kernel function $\cot \left( t/2\right) :$
singular ($2/t)$, and regular. Fast summation with the singular kernel can
be performed using the MLFMM, while the regular part can be truncated and
factorized for all points in the respective neighborhoods at level 2. The
MLFMM for kernel $1/t$ has been well studied, optimized, and error bounds
are established \cite{GDB03}. We just modified our software in \cite{GDB03}
with a circular data structure.

\subsection{Truncation Errors:} 
Consider the truncation error for sum (\ref{10})
when $q$ first terms of the expansion are used to approximate the infinite
sum. We note that the Bernoulli numbers satisfy \cite{AS64} 
\begin{equation}
\left\vert B_{2m}\right\vert <\frac{2\left( 2m\right) !}{\left( 2\pi \right)
^{2m}}\frac{1}{1-2^{1-2m}},\quad m=1,2,...  \label{12}
\end{equation}%
This shows that the truncated part of the series can be majorated by the
geometric progression. Using $\left\vert t\right\vert \leqslant \pi ,$ the
error bound is%
\begin{equation}
\left\vert \epsilon _{q}^{(1)}\right\vert \!=\left\vert
\!\sum_{m=q+1}^{\infty }\frac{\left\vert B_{2m}\right\vert t^{2m-1}}{\left(
2m\right) !}\right\vert <\frac{2^{1-2q}}{3\pi \left( 1-2^{-2q-1}\right) }%
=\epsilon _{q}.  \label{p14}
\end{equation}%
The truncation error, $\epsilon _{q}^{(2)},$ in (\ref{11}) can be bounded
similarly, using $\left\vert t\right\vert \leqslant \pi /2$: $\left\vert
\epsilon _{q}^{(2)}\right\vert <\!2\epsilon _{q}.$ The estimates are for a
single source-evaluation pair. For $3N/4$ sources located in the
neighborhood of point $y_{j}$, and $N/4$ outside it, the total maximum
absolute error in computation of $g_{j}$ (with $\left\vert f_{k}\right\vert
\leqslant 1)$ is 
\begin{equation}
\!\!\!\!\!\!\left\vert \epsilon _{q}^{reg}\right\vert \!=\!\left\vert
\sum_{k=0}^{N-1}\left( K_{jk}-K_{jk}^{(q)}\right) f_{k}\!\right\vert <\frac{2%
}{N}\left( \frac{3N\epsilon _{q}}{4}+\frac{2N\epsilon _{q}}{4}\right) =\frac{%
5\epsilon _{q}}{2}.  \notag
\end{equation}%
The error of the MLFMM with kernel $1/t$ has been estimated in \cite{GDB03}.
This result applied to our case in a computational domain of size $3\pi /2, $
and for $3N/4$ sources is 
\begin{equation}
\left\vert \epsilon _{p}^{MLFMM}\right\vert <\frac{5}{\pi }\cdot 2^{l_{\max
}}3^{-p},  \label{17}
\end{equation}%
where $l_{\max }$ is the maximum level of space subdivision and $p$ is the
truncation number used in the MLFMM. Therefore the total truncation error
for the present method can be estimated as 
\begin{equation}
\epsilon \lesssim \frac{5}{3\pi }\left( 4^{-q}+2^{l_{\max }}3^{1-p}\right) .
\label{18}
\end{equation}%
By requiring that the error in the singular and regular parts be the same in
(\ref{18}) we relate \thinspace $p$ and $q$, and obtain a combined bound as 
\begin{equation}
q\gtrsim \frac{1}{2}\log _{2}\frac{3\pi }{10\epsilon },\quad p\gtrsim \log
_{3}\frac{3\pi }{10\epsilon }+\frac{l_{\max }}{\log _{2}3}+1.  \label{18.2}
\end{equation}

\subsection{Complexity and Optimizations: }
For fast computations of the products
involving $q$-truncated sums (\ref{10}) and (\ref{11}) we factorize powers $%
t^{2m-1}$ using the Newton binomial expansion near the center, $x_{c}^{(n)}$%
, of the segment $P_{n}$ containing $y_{j}.$ Substituting this expansion
into Eqs (\ref{10}) and (\ref{11}) and further into Eq. (\ref{p9}), we
obtain after changing the order of summation: 
\begin{gather}
h\left( y_{j}\right) =2\sum_{x_{k}\in \Omega _{1}(P_{n})}\frac{f_{k}}{y_{j}-%
\widetilde{x}_{k}}-\sum_{l=0}^{2q-1}\frac{d_{l}}{l!}\left(
y_{j}-x_{c}^{(n)}\right) ^{l},  \label{p20} \\
d_{l}=\sum_{m=\left[ l/2\right] +1}^{q}\frac{\left\vert B_{2m}\right\vert }{m%
}\left[ \alpha _{2m-l-1}^{(1)}+\left( 2^{2m}-1\right) \alpha _{2m-l-1}^{(2)}%
\right] ,  \notag \\
\alpha _{l}^{(s)}=\frac{1}{l!}\sum_{x_{k}\in \Omega _{s}(P_{n})}f_{k}\left(
x_{c}^{(n)}-\widetilde{x}_{k}\right) ^{l},\quad s=1,2.  \notag
\end{gather}%
For $2q\ll \min (N,M)$ the second sum in Eq. (\ref{p20}) is $O\left(
2q\left( N+M\right) \right) .$ For the complexity of the MLFMM in the
present case we have \cite{GDB03} 
\begin{equation}
C_{FMM}=O\left( p\left( N+M\right) +\frac{3NM}{2^{l_{\max }}}+\frac{6P}{%
2^{-l_{\max }}}\right) ,  \label{23}
\end{equation}%
where $P(p)$ is the cost of a single translation for truncation $p$ .

The total cost of the fast Fourier interpolation algorithm (FFIA) can be
estimated, and minimized by selection of $l_{\max }$ for given error (\ref%
{18}). A simplified estimate assuming that $p$ and $q$ change slower than $%
2^{l_{\max }},$ yields 
\begin{equation}
l_{\max }^{(opt)}\sim \frac{1}{2}\log _{2}\frac{NM}{2P},  \label{25}
\end{equation}%
and the total complexity of the optimized algorithm will be 
\begin{equation}
C_{FFIA}^{(opt)}=O\left( \left( N+M\right) \left( p+2q\right) +6\left[
2NMP(p)\right] ^{1/2}\right) .  \label{26}
\end{equation}%
For $N\sim M,$ this yields $C_{FFIA}^{(opt)}=O\left( N\left( \log N+\log
\epsilon ^{-1}\right) \right) $.%

\section{Numerical study}

The FFIA was implemented and tested on a 933 MHz Pentium III\ Xeon PC with
1GB RAM. in complex double precision. Tests were performed for $M$ $=N=2^{L}$
varying in $2^{3}-2^{20}$ for various $l_{\max },p,$ and $q$. In all cases, $%
y_{j},$ $j=0,...,M-1$ were randomly distributed in $\left[ 0,2\pi \right) $
for the transform $\left\{ f_{k}\right\} \rightarrow \left\{ g_{j}\right\} $
and were within 10\% perturbation near the centers $\left\{ x_{k}\right\} $
for the transform $\left\{ g_{j}\right\} \rightarrow \left\{ f_{k}\right\} .$
We found that the errors for the latter depend substantially on the
distribution of $\left\{ y_{j}\right\} ,$ while for $\left\{ f_{k}\right\}
\rightarrow \left\{ g_{j}\right\} $ there was no such dependence. We thus
focused on the transform $\left\{ f_{k}\right\} \rightarrow \left\{
g_{j}\right\} ,$ since the error analysis for the other case should be
combined with a more rigorous study of non-uniformity errors and errors in
computation of coefficients (\ref{p8}).

\begin{figure}[htb]
	\includegraphics[width=0.9\textwidth,trim=0in 0in 1in 0in]{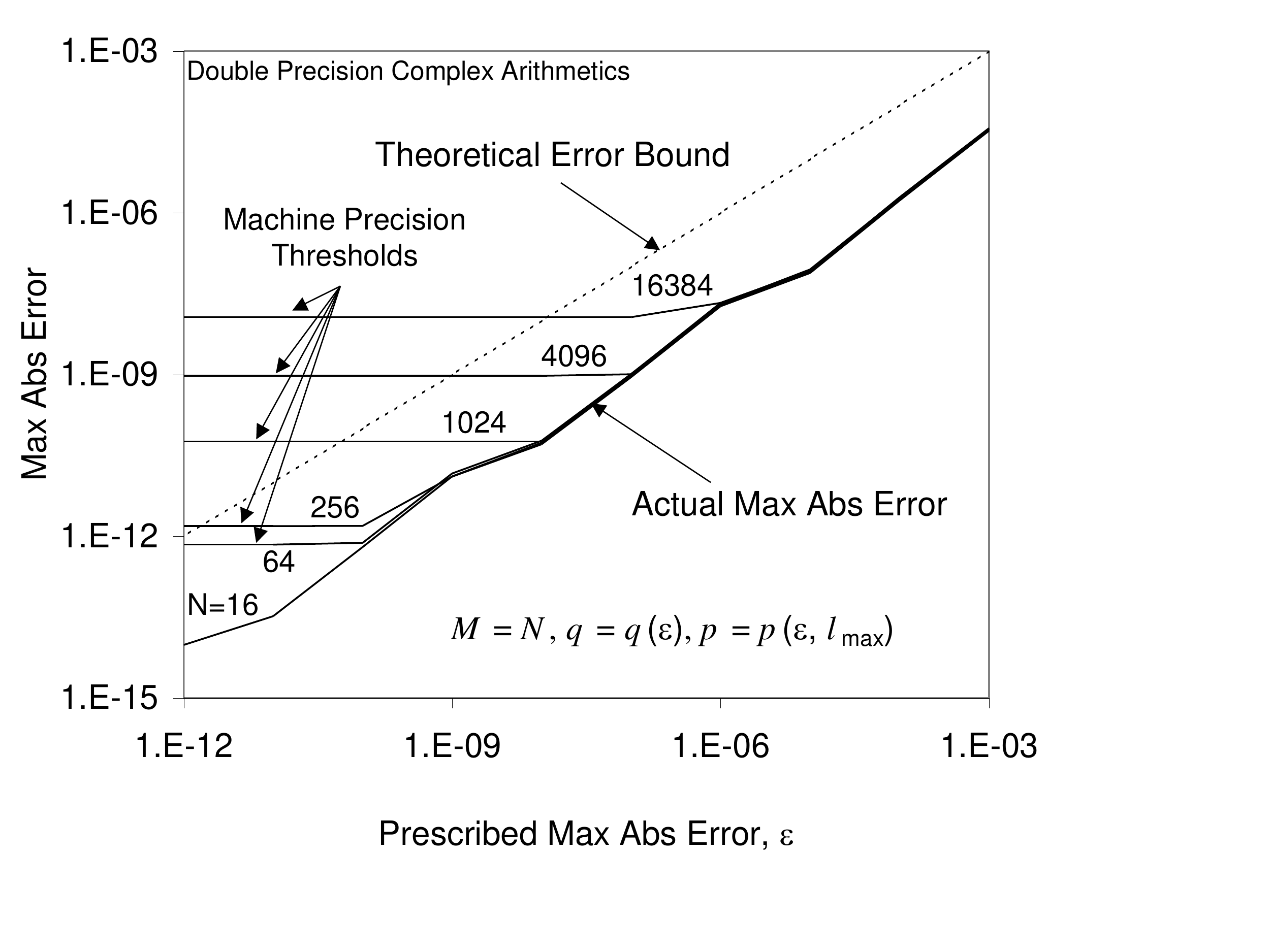}
	\caption{Dependences of $\protect\epsilon _{a}$, for the FFIA solution on
		the prescribed $\protect\epsilon $ for various $N$, indicated near the
		curves.  \label{Fig.comp1}}
\end{figure}

\subsection{Error analysis}

We performed several error tests to ensure that the algorithm worked
properly and results were consistent with theory. First, we took the
analytical solution $f\left( x\right) \equiv 1,.$ and varied $l_{\max }$
from $3$ to $L=\log _{2}N$, and prescribed $\epsilon =10^{-12},...,10^{-3};$ 
$q$ and $p$ were found from (\ref{18.2}). Results of the tests are in Fig. %
\ref{Fig.comp1}. The actual error was measured as $\epsilon
_{a}=\max_{j}\left\vert g_{j}-1\right\vert $ and depends on $N,l_{\max
},\epsilon .$ We observed that $\epsilon _{a}$ just slightly depends on $%
l_{\max }$, consistent with (\ref{18}). Another theoretical prediction
provided by Eq. (\ref{18}) is that the error should not depend on $N$. Fig. %
\ref{Fig.comp1} shows this for $\epsilon \gg \epsilon _{th}(N),$ where $%
\epsilon _{th}(N)$ is a threshold error, which we refer as \textquotedblleft
machine precision error\textquotedblright ; $\epsilon _{a}$ was much smaller
than $\epsilon $.
\begin{figure}[htb]
	\includegraphics[width=0.9\textwidth,trim=0in 0in 1in 0in]{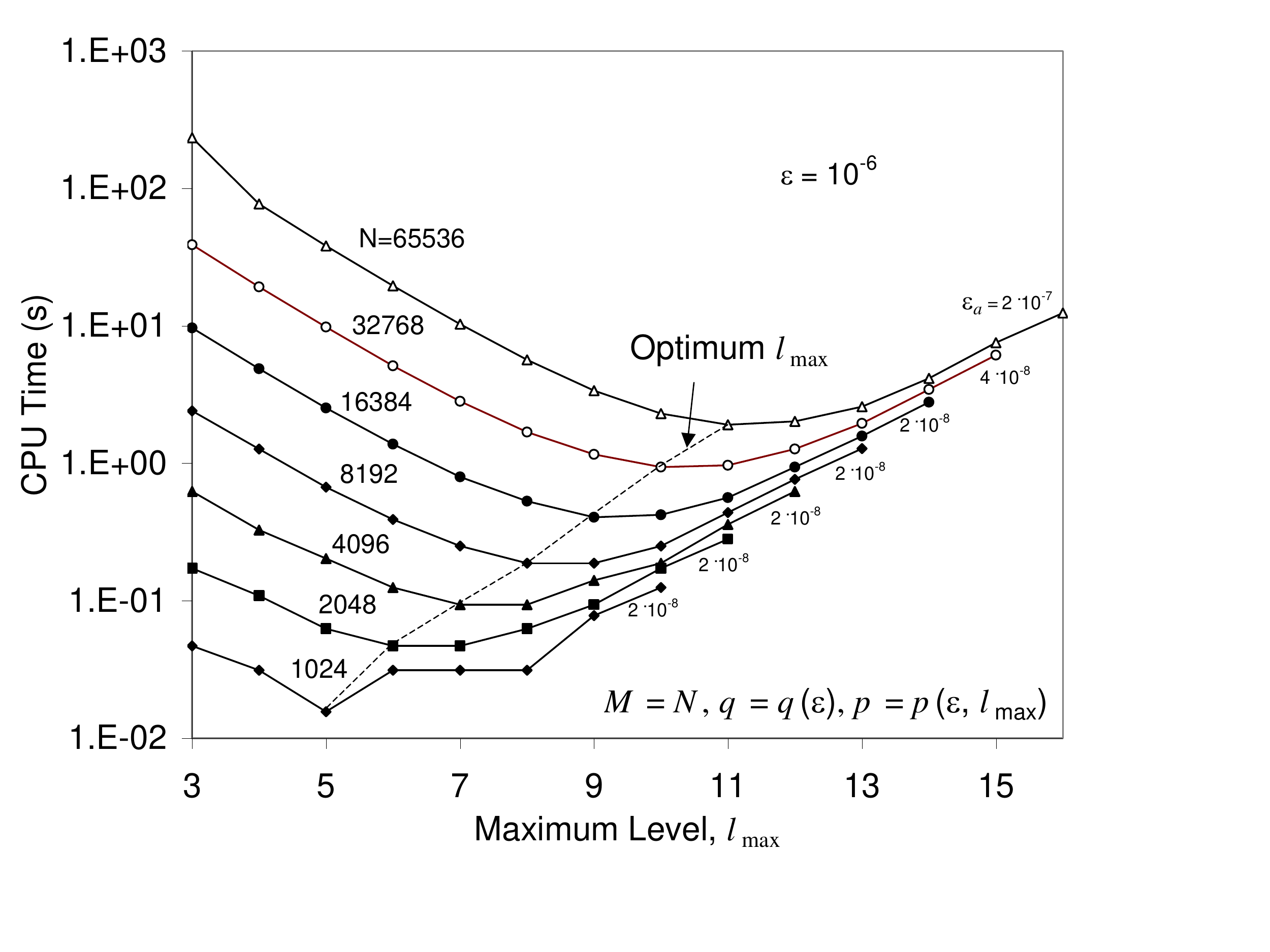}
	\caption{Dependence of the CPU time on $l_{\max }$ in the MLFMM for prescribed $\protect\epsilon =10^{-6}$ for different $N$. The small numbers near the right end of the curves show $\protect\epsilon _{a}$  \label{Fig.comp4}}
\end{figure}

To validate that $\epsilon _{th}$ is really related to the machine precision
arithmetic, we performed a test, when $\left\{ g_{j}\right\} $ was computed
by the \textquotedblleft exact\textquotedblright\ straightforward method,
i.e. by multiplication of matrix $\left\{ K_{jk}\right\} $ computed using
Eq. (\ref{p4}) by vector $\left\{ f_{k}\right\} $. The results of this test
yield dependence $\epsilon (N),$ which is consistent with $\epsilon _{th}(N)$
obtained in experiments using the \textquotedblleft
approximate\textquotedblright\ fast method. We also compared the maximum
absolute difference between the results obtained using the straightforward
method and the FFIA for constant and random $f_{k}\in \left[ 0,1\right] $.
In both cases for $\epsilon =10^{-12}$ this error was smaller than $\epsilon
_{th}(N).$

We can thus conclude that the theoretical error or machine precision bound
the FFIA\ error. It makes no sense to specify an error below the precision,
and this provides a result on the asymptotic complexity of the MLFMM.
Indeed, estimation $C_{FFIA}^{(opt)}=O(N\log N)$ does not take into account
the growth of $\epsilon _{th}(N)$ with $N$ due to roundoff errors, while
assumes that computations are done with arbitrary precision. However, this
happens only for $\epsilon \gg \epsilon _{th}(N)$. If the computational task
is formulated as \textquotedblleft solve the problem with machine
precision\textquotedblright\ then the complexity will be $O(N).$

\subsection{Optimization:} We conducted optimization tests in two basic settings
for error control. The first case we considered is minimization of the CPU
time for fixed prescribed error $\epsilon \geqslant \epsilon _{th}(N)$. The
second case is minimization of the CPU time for computations with maximum
available machine precision. In the first case the range of $N$ was limited
by condition $\epsilon \geqslant \epsilon _{th}(N)$ and in the second case $%
N $ was limited by machine resources available.

\begin{figure}[htb]
	\includegraphics[width=0.9\textwidth,trim=0in 0in 1in 0in]{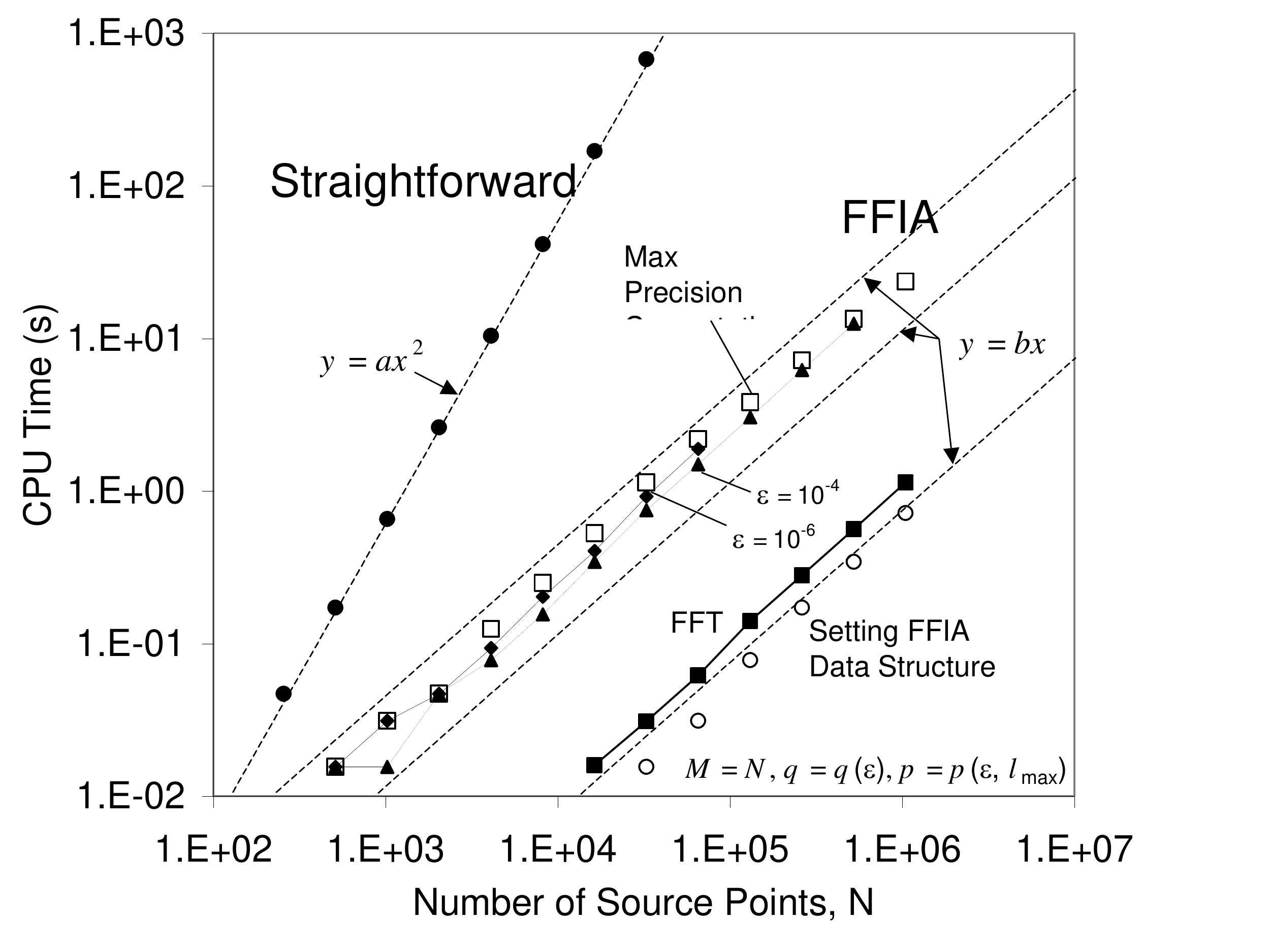}
	\caption{Dependences of the CPU time
		on $N,$ for the straightforward method (dark circles), and present method
		(FFIA) with prescribed accuracy (dark diamonds and triangles) and with the
		machine precision (light squares). The light circles show CPU time for
		setting the data structure. The dark squares show the CPU time required for
		FFT on the same machine.  \label{Fig.comp5}}
\end{figure}

Fig. \ref{Fig.comp4} illustrates dependence of the CPU time on $l_{\max }$
at fixed $\epsilon =10^{-6}$ and different $N$. Truncation numbers $q$ and $%
p $ were selected using (\ref{18.2}). The actual error of computations $%
\epsilon _{a}$ was estimated by comparison with analytical solution, $%
f\equiv 1$ and shown on Fig. \ref{Fig.comp4}.

\begin{figure}[htb]
	\includegraphics[width=0.9\textwidth,trim=0in 0in 1in 0in]{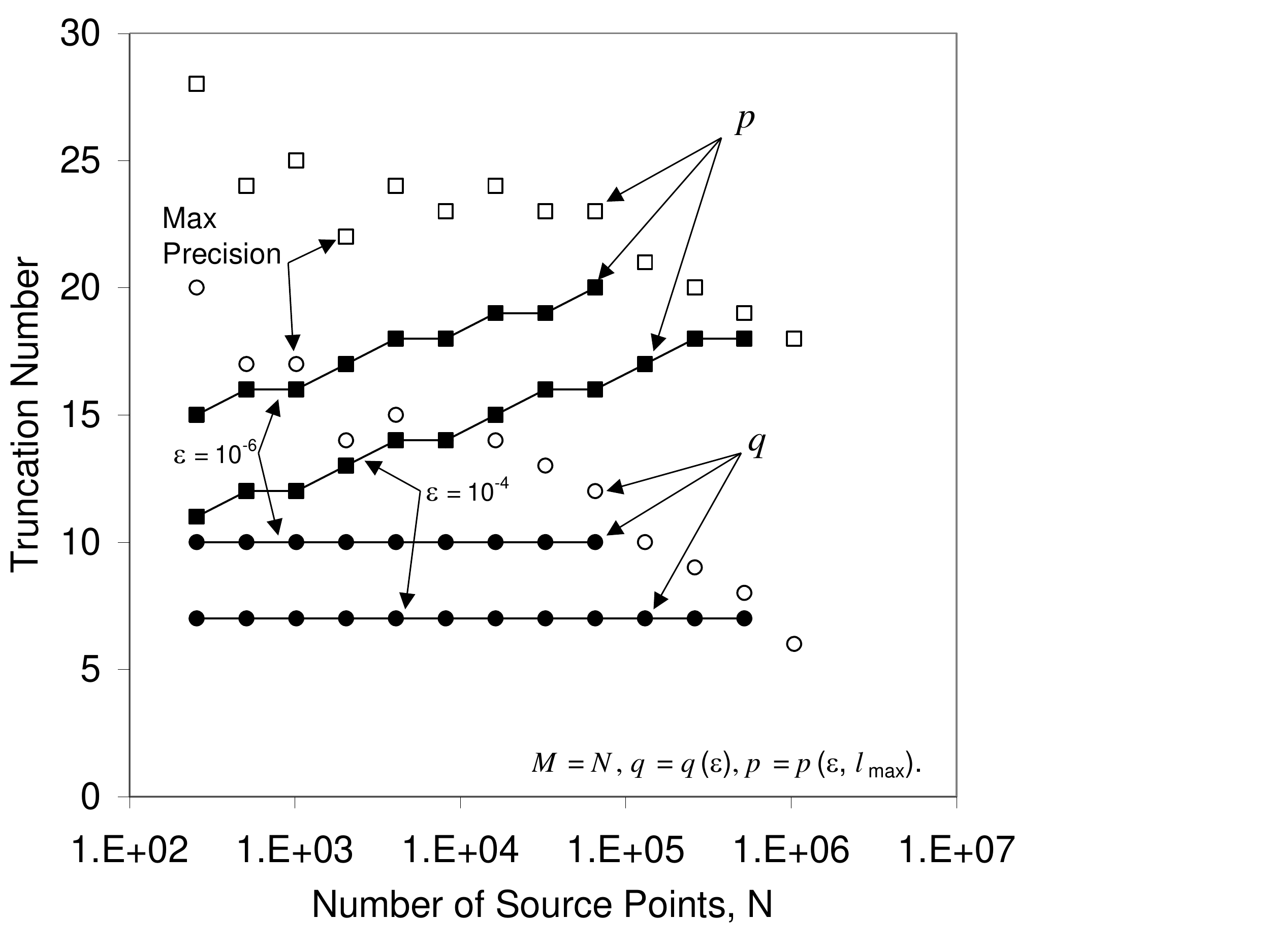}
	\caption{Truncation numbers, $q$ (circles) and $p$ (squares), vs $N$. Dark circles and squares show
		dependence of computation with fixed $\protect\epsilon $. Light circles and squares show the threshold dependences 
		$q_{th}\left(N\right),$ $p_{th}\left( N\right)$ Obtained using Eq. (\ref{18.2}) for $\epsilon =\epsilon _{th}\left( N\right) $.\label{Fig.comp6}}
\end{figure}

All curves
on Fig. \ref{Fig.comp4} clearly show that there exist an optimum maximum
level of the space subdivision, $l_{\max }^{(opt)},$ for the MLFMM. The CPU
time for computations with $l_{\max }=l_{\max }^{(opt)}$ and non-optimized $%
l_{\max }$ can differ substantially. An empirical correlation is $l_{\max
}^{(opt)}=\log _{2}N-l_{\ast }(\epsilon ),$ where $l_{\ast }(\epsilon )$
depends on the algorithm implementation (in our case $l_{\ast }(\epsilon )=5$%
, see Fig. \ref{Fig.comp4}).

In the second case we imposed experimentally measured machine precision
threshold $\epsilon _{th}(N)$ as the prescribed accuracy $\epsilon .$ In
this case we found that the same dependence $l_{\max }^{(opt)}(N)$ describes 
$l_{\max }^{(opt)}\left( L\right) $ well. We also noted that the actual
error of computations $\epsilon _{a}$ was equal to $\epsilon _{th}(N)$ for $%
\epsilon =\epsilon _{th}(N).$

\subsection{Speed of Computations: }To evaluate the algorithm we measured CPU
time for the straightforward method, and the optimized FFIA for computations
with fixed $\epsilon \geqslant \epsilon _{th}(N)$ and machine precision, $%
\epsilon =\epsilon _{th}(N).$ It is noticeable that the MLFMM procedures
consist of two steps: setting the data structure, the step which needs to be
performed when the target data set $\left\{ y_{j}\right\} $ changes, and
evaluation, which is executed each time the input $\left\{ f_{k}\right\} $
changes. Fig. \ref{Fig.comp5} shows that the computational cost for setting
the data structure grows as $O(N)$ and is much smaller than the cost of the
evaluation step. In case the interpolation is performed for INUFFT, an
execution of the IFFT is required for the transform $\left\{ c_{n}\right\}
\rightarrow \left\{ f_{k}\right\} $. Fig. \ref{Fig.comp5} shows that the
cost of this step is low compared to the FFIA evaluation step and is
slightly above that of setting of the data structure.

Fig. \ref{Fig.comp5} also shows that the FFIA far outperforms the
straightforward method even for small data sets ($N\sim 10^{2}$) and the
difference gains orders of magnitude for larger $N.$ This is clear, since
the straightforward method is scaled as $O(N^{2}),$ while the FFIA for
maximum precision is scaled as $O(N)$. It is interesting to see that for
fixed prescribed error $\epsilon ,$ such as $\epsilon =10^{-6}$ and $%
\epsilon =10^{-4}$ shown in the figure, the CPU time required for
computations deviates from the linear dependence $O(N)$ at larger $N.$ This
is consistent with the estimation of the algorithm complexity $O(N\log N)$
valid for for $\epsilon \gg \epsilon _{th}(N).$ In contrast to this behavior
the curve corresponding to computations with the maximum precision
available, $\epsilon =\epsilon _{th}(N),$ shows deviation from the linear
dependence in the\emph{\ opposite} direction at larger $N$. The explanation
of this effect comes from the error analysis provided above and also clear
from Fig. \ref{Fig.comp6}, which illustrates the dependences of the
truncation numbers on $N.$ While $q$ for fixed $\epsilon $ stays constant
and $p$ increases $a+b\log N$ a different behavior of the bounding values $%
p_{th}$ and $q_{th}$ for maximum precision computations occurs. Both the
truncation numbers, $p_{th}$ and $q_{th}$, are bounded by a global constant.
Moreover, they have a tendency to decay at larger $N$ and $p_{th},q_{th}\sim
a-b\log N,$ where $a$ and $b$ some constants.%

\section{Conclusions}

We developed and tested a fast algorithm for matrix-vector multiplication
with a kernel appropriate for the forward and inverse interpolation problems
of band limited functions, which can be used for the forward and inverse
Fourier Transforms for non-equispaced data. The algorithm far outperforms
straightforward methods. Our results stay within the specified error bounds
or maximum available machine accuracy. The algorithm asymptotically scales
as $O(N)$ for computations with fixed maximum machine precision and as $%
O(N\log N+N\log \epsilon ^{-1})$ otherwise. Furthermore, the current
algorithm is parallelizable.

\end{document}